# Numerical method based on Galerkin approximation for the fractional advection-dispersion equation

Harendra Singh, Manas Ranjan Sahoo, Om Prakash Singh

Department of Mathematical Sciences, Indian Institute of Technology, Banaras Hindu University, Varanasi 221005, India

**Abstract:** We use a concept of weak asymptotic solution for homogeneous as well as non-homogeneous fractional advection dispersion type equations. Using Legendre scaling functions as basis, a numerical method based on Galerkin approximation is proposed. This leads to a system of fractional ordinary differential equations whose solutions in turn give approximate solution for the advection-dispersion equations of fractional order. Under certain assumptions on the approximate solutions, it is shown that this sequence of approximate solutions forms a weak asymptotic solution. Numerical examples are given to show the effectiveness of the proposed method.

## 1 Introduction

We consider the following advection-dispersion equation

$$\frac{\partial u}{\partial t} = -v\frac{\partial u}{\partial x} + k\frac{\partial^2 u}{\partial x^2}, \tag{1}$$

where $u$ is the solute concentration in a fluid, the positive constants $v, k$ represent the average fluid velocity and dispersion coefficient respectively, $x$ is spatial domain and $t$ is the time variable. This system is used to describe diffusion as well as transport phenomenon such as concentration etc. A rich theory exists in the literature for existence and uniqueness of the solution and its continuous dependence on the initial data for this equation when derivatives are of integer order. If integer orders are replaced by fractional order, then not much is known in this direction. Fractional derivatives are non-local in nature whereas integer orders are local. So it is physically very important to study such systems because as it is a well- known fact that variation of some physical quantities depends on the past.

---

E-mail addresses: opsingh.apm@iitbhu.ac.in; harendrasingh.rs.apm12@iitbhu.ac.in; manastifr@gmail.com;

There exist a vast literature for theory and applications of fractional differential equations in areas such as finance [1-3], physics [4-7], control theory [8] and hydrology [9-13].Several papers have been written [14-16] to show the equivalence between the transport equations using fractional order derivatives and some heavy-tailed motions.

There exist several methods to solve fractional advection-dispersion equation such as implicit and explicit difference method [17-20], Green function [10], variable transformation [21], Adomian's decomposition method (ADM) [22] and optimal homotopy asymptotic method [23].

Like the integer order advection-dispersion equation, there does not exist any analytical method for the existence and uniqueness of the solution. Only approximate solutions may be obtained by the existing methods.

We study the non-homogeneous advection dispersion equation of fractional order, namely,

$$\frac{\partial^\alpha u(x,t)}{\partial t^\alpha} = -v \frac{\partial^\beta u(x,t)}{\partial x^\beta} + k \frac{\partial^\gamma u(x,t)}{\partial x^\gamma} + f(x,t),$$

$$t, x > 0; 0 < \alpha, \beta \leq 1; 1 < \gamma \leq 2, \tag{2}$$

Subject to the initial conditions:

$$u(x,0) = g(x). \tag{3}$$

In this method we first get finite dimensional approximate solution taking finite dimensional basis in x-direction and varying coefficients as functions of t, which in turn leads a system of ODE having constant coefficients whose solution is known [26]. In this way we get the coefficients, so the approximate solution.

The present paper is organised as follows. In section-2, we describe basic preliminaries for fractional differentiation, Legendre scaling functions and introduce concept of solution. In section-3, we give description of the algorithm for the construction of approximate solutions. In section-4, we show existence of a weak asymptotic solution under certain conditions. In section-5, we present numerical experiment to show the effectiveness of the proposed method.

## 2 Preliminaries

There are several definitions of fractional order derivatives and integrals. These are not necessarily equivalent. In this paper, the fractional order differentiations and integrations are in the sense of Caputo and Riemann-Liouville sense respectively.

**Definition 2.1** The Riemann-Liouville fractional order integral operator is given by [24]

$$I^\alpha f(x) = \frac{1}{\overline{|\alpha}} \int_0^x (x-t)^{\alpha-1} f(t) dt \quad \alpha > 0, x > 0,$$

$$I^0 f(x) = f(x).$$

For the fractional Riemann-Liouville integration

$$I^\alpha x^k = \frac{\Gamma(k+1)}{\Gamma(k+1+\alpha)} x^{k+\alpha}$$

**Definition 2.2** The Caputo fractional derivative of order $\beta$ are defined as

$$D^\beta f(x) = I^{m-\beta} D^m f(x) = \frac{1}{\Gamma(m-\beta)} \int_0^x (x-t)^{m-\beta-1} \frac{d^m}{dt^m} f(t) dt,$$

$m-1 < \beta < m, x > 0$.

For the Caputo derivative, we have [25]

$D^\beta A = 0$ (A is a constant),

$$D^\beta x^k = \frac{\Gamma(k+1)}{\Gamma(k+1-\beta)} x^{k-\beta} \text{ for } k \in \mathbb{N}_0 \text{ and } k \geq \lceil \beta \rceil \text{ or } k \notin \mathbb{N} \text{ and } k > \lfloor \beta \rfloor,$$

$= 0$ for $k \in \mathbb{N}_0$ and $k < \lceil \beta \rceil$

Where $\lceil \beta \rceil$ and $\lfloor \beta \rfloor$ are the ceiling and floor functions respectively, while $\mathbb{N} = \{1, 2, 3, \ldots\}$ and $\mathbb{N}_0 = \{0, 1, 2, \ldots\}$.

The Caputo fractional differentiation is a linear operator similar to the integer order differentiation.

The Legendre scaling functions $\{\phi_i(x)\}$ are defined by

$$\phi_i(x) = \begin{cases} \sqrt{(2i+1)} P_i(2x-1), & \text{for } 0 \leq x < 1. \\ 0, & \text{otherwise} \end{cases}$$

Where $P_i(x)$ is Legendre polynomials of order $i$ which are orthogonal on the interval $[-1, 1]$ with respect to the weight function $w(x) = 1$. Legendre scaling functions are constructed normalizing the shifted Legendre polynomials. So the collections $\{\phi_i(x)\}$ form an orthonormal basis for $L^2[0,1]$. A function $f \in L^2([0,1] \times [0, \infty])$, with bounded second derivative $\left| \frac{\partial^2 f(x,t)}{\partial x^2} \right| \leq M$, expanded as infinite sum of Legendre scaling functions for fixed value of $t$, converges uniformly to the function $f(x,t)$, independent of the variables $t$ and $x$ lying in a compact set.

$$f(x,t) = \lim_{n \to \infty} \sum_{i=0}^{n} c_i(t) \phi_i(x), \tag{4}$$

where $c_i(t) = \langle f(t,x), \phi_i(x) \rangle$, and $\langle .,. \rangle$ is standard inner product on $L^2[0,1]$.

If the series is truncated at $n = m$, then we have

$$f \cong \sum_{i=0}^{m} c_i \phi_i = c^T(t) \phi(x), \tag{5}$$

where, $C(t)$ and $\phi(x)$ are $(m+1) \times 1$ matrices given by

$$c(t) = [c_0(t), c_1(t), \ldots, c_m(t)]^T \text{ and } \phi(x) = [\phi_0(x), \phi_1(x), \ldots, \phi_m(x)]^T.$$

The Legendre scaling functions of degree $i$ are given by

$$\phi_i(x) = (2i+1)^{\frac{1}{2}} \sum_{k=0}^{i} (-1)^{i+k} \frac{(i+k)!}{(i-k)!} \frac{x^k}{(k!)^2}, \tag{6}$$

where

$$D^\beta \phi_i(x) = 0, \qquad i = 0, 1, \ldots, \lceil \beta \rceil - 1, \quad \beta > 0. \tag{7}$$

First of all we define the concept of weak asymptotic solution introduced by Panov and Shelkovich [27], for its application, see, [28] and [29].

**Definition:** A sequence of function $u_n$ is said to be weak asymptotic solution if the following holds

$$\lim_{n \to \infty} \int_0^1 \left( \frac{\partial^\alpha u_n(x,t)}{\partial t^\alpha} + v \frac{\partial^\beta u_n(x,t)}{\partial x^\beta} - k \frac{\partial^\gamma u_n(x,t)}{\partial x^\gamma} - f_n(x,t) \right) \phi(x) dx = 0, \text{ uniformly in } t \text{ lying in the}$$

compact subset of $(0, \infty)$ and $\lim_{n \to \infty} \int_0^1 (u_n(x,0) - u_0(x)) \phi(x) dx = 0, \text{ as } n \to \infty.$

where $\phi$ is a $C^\infty$ function having compact support in $(-\infty, \infty)$.

## 3. Method of solution

Now, we consider the following non-homogeneous fractional advection-dispersion equation,

$$\frac{\partial^\alpha u(x,t)}{\partial t^\alpha} = -v \frac{\partial^\beta u(x,t)}{\partial x^\beta} + k \frac{\partial^\gamma u(x,t)}{\partial x^\gamma} + f(x,t), \tag{8}$$

$t, x > 0; 0 < \alpha, \beta \leq 1; 1 < \gamma \leq 2,$

with initial condition $u(x,0) = g(x).$ (9)

Using (5), we approximate $u(x,t)$ by

$$u_n(x,t) = \sum_{i=0}^{n} c_i(t)\phi_i(x) = c^T(t)\phi(x), \qquad (10)$$

and taking Caputo fractional derivatives of order $\alpha$ and $\beta$ of Eq. (10) with respect to t and x respectively, we get

$$\frac{\partial^\alpha u_n(x,t)}{\partial t^\alpha} = \sum_{i=0}^{n} c_i^\alpha(t)\phi_i(x), \qquad (11)$$

$$\frac{\partial^\beta u_n(x,t)}{\partial x^\beta} = \sum_{i=0}^{n} c_i(t)\phi_i^\beta(x). \qquad (12)$$

Again truncating $\frac{\partial^\beta u_n(x,t)}{\partial x^\beta}$ using first n terms of Legendre scaling expansion, in equation Eq. (12), we get

$$\frac{\partial^\beta u_n(x,t)}{\partial x^\beta} = \sum_{i=0}^{n} c_i(t)\phi_i^\beta(x) = \sum_{j=0}^{n} d_j(t)\phi_j(x), \qquad (13)$$

The coefficients are given by

$$d_j(t) = \int_0^1 \sum_{i=0}^{n} c_i(t)\phi_i^\beta(x)\phi_j(x)dx = \sum_{i=0}^{n} c_i(t)a_{ij}, \qquad (14)$$

where

$$a_{ij} = \int_0^1 \phi_i^\beta(x)\phi_j(x)dx, \quad i = 0,1,....,n, \text{ and } j = 0,1,....,n. \qquad (15)$$

Similarly taking derivative of order $\gamma$ of Eq. (10) with respect to x and truncating first n terms of the expansion, we get

$$\frac{\partial^\gamma u_n(x,t)}{\partial x^\gamma} = \sum_{i=0}^{n} c_i(t)\phi_i^\gamma(x) = \sum_{j=0}^{n} e_j(t)\phi_j(x), \qquad (16)$$

The coefficients are given by

$$e_j(t) = \int_0^1 \sum_{i=0}^{n} c_i(t)\phi_i^\gamma(x)\phi_j(x)dx = \sum_{i=0}^{n} c_i(t)b_{ij}, \qquad (17)$$

where

$$b_{ij} = \int_0^1 \phi_i^\gamma(x)\phi_j(x)dx, \ i = 0,1,....,n, \text{ and } j = 0,1,....,n. \tag{18}$$

The non-homogeneous term $f(x,t)$ is also approximated as,

$$f_n(x,t) = \sum_{i=0}^n r_i(t)\phi_i(x). \tag{19}$$

Now putting the approximate values of $\dfrac{\partial^\alpha u_n(x,t)}{\partial t^\alpha}, \dfrac{\partial^\beta u_n(x,t)}{\partial x^\beta}, \dfrac{\partial^\gamma u_n(x,t)}{\partial x^\gamma}$ and $f_n(x,t)$ from Eq. (11), (12), (16) and (19) in Eq. (8), we get,

$$\sum_{i=0}^n c_i^\alpha(t)\phi_i(x) = -v\sum_{j=0}^n d_j(t)\phi_j(x) + k\sum_{j=0}^n e_j(t)\phi_j(x) + \sum_{i=0}^n r_i(t)\phi_i(x), \tag{20}$$

Substituting the value of $d_j(t)$ and $e_j(t)$ from Eq. (14) and Eq. (17) in Eq. (20) we get,

$$\sum_{i=0}^n c_i^\alpha(t)\phi_i(x) = -v\sum_{j=0}^n\sum_{i=0}^n c_i(t)a_{ij}\phi_j(x) + k\sum_{j=0}^n\sum_{i=0}^n c_i(t)b_{ij}\phi_j(x) + \sum_{i=0}^n r_i(t)\phi_i(x),$$

Comparing the coefficients of scaling functions in above equation yields the following system of differential equations,

$$c^\alpha(t) = Ac(t) + r(t), \tag{21}$$

Where $A = (kb_{ij} - va_{ij})$, $i = 0,1,....,n$, $j = 0,1,....,n$. is a $(n+1)\times(n+1)$ matrix, $c(t) = (c_0(t), c_1(t),......,c_n(t))^T$, $c^\alpha(t) = (c_0^\alpha(t), c_1^\alpha(t),......,c_n^\alpha(t))^T$ and $r(t) = (r_0(t), r_1(t),......,r_n(t))^T$.

The solutions of the above differential equation for integer order as well as fractional order are discussed below.

**Case 1.** When $\alpha = 1$, then the above system becomes

$$c'(t) = Ac(t) + r(t), \tag{22}$$

whose solution is given by

$$c(t) = e^{At}c(0) + e^{At}\int_0^t e^{-As}r(s)ds, \tag{23}$$

where $c(0) = (c_0(0), c_1(0),......,c_n(0))^T$.

The coefficients $c_i(0), \quad i = 0,1,....,n,$ can be determined from the approximation of initial value $u(x,0)$.

$$u(x,0) = g(x) \cong \sum_{i=0}^{n} c_i(0)\phi_i(x), \tag{24}$$

this implies $c_i(0) = \int_0^1 g(x)\phi_i(x)dx. \tag{25}$

Substituting the value of $c(t)$ from Eq. (23) in Eq. (10), approximate solutions for non-homogeneous fractional advection dispersion can be obtained.

If system of differential equation is homogeneous then Eq. (23) will become

$$c(t) = e^{At}c(0) \tag{26}$$

Substituting the value of $c(t)$ from Eq. (26) in Eq. (10), approximate solutions for homogeneous fractional advection dispersion can be obtained.

**Case 2.** When $0 < \alpha < 1$

$$c(t) = e_\alpha^{At}c(0) + \int_0^t e_\alpha^{A(t-s)} r(s)ds \tag{27}$$

Where $\alpha$ – exponential function, is defined by

$$e_\alpha^{A(t-a)} = (t-a)^{\alpha-1} \sum_{i=0}^{\infty} A^i \frac{(t-a)^{i\alpha}}{|\{(i+1)\alpha\}|}. \quad [26] \tag{28}$$

Substituting the value of $c(t)$ from Eq. (28) in Eq. (10), approximate solutions for non-homogeneous fractional advection dispersion can be obtained.

If system of differential equation is homogeneous then Eq. (27) will become

$$c(t) = e_\alpha^{At}c(0) \tag{29}$$

Substituting the value of $c(t)$ from Eq. (29) in Eq. (10), approximate solutions for homogeneous fractional advection dispersion can be obtained.

## 4. Convergence Analysis

**Lemma:** Let $f$ be a $L^2$ - function whose second derivative is bounded i.e. $|D^2 f| \leq K$ and $S_n(f)$ be the n<sup>th</sup> Fourier sum using Legendre scaling function, then we have the following estimates

$$\|f(x) - S_n f(x)\|^2_{L^2} \leq \frac{K^2}{256} F_3(-\frac{3}{2} + n) \tag{30}$$

where,

$$\|f(x)\|_{L^2} = \left( \int_0^1 |f(x)|^2 \, dx \right)^{\frac{1}{2}}.$$

Where $F_n(z)$ is the PolyGamma function defined by,

$$F_n(z) = (-1)^{n+1} \lfloor n \sum_{k=0}^{\infty} \frac{1}{(z+k)^{n+1}}.$$

**Proof:** Let $f(x,t) = \sum_{i=0}^{\infty} c_i(t) \phi_i(x)$. Truncating it to level $n-1$, we get $S_n f(x) = \sum_{i=0}^{n-1} c_i(t) \phi_i(x)$, thus,

$$f(x) - S_n f(x) = \sum_{i=n}^{\infty} c_i(t) \phi_i(x).$$

(31)

$$\|f(x) - S_n f(x)\|^2_{L^2} = \int_0^1 (f(x) - S_n f(x))^2 \, dx$$

$$= \int_0^1 \left( \sum_{i=n}^{\infty} c_i(t) \phi_i(x) \right)^2 dx$$

$$= \sum_{i=n}^{\infty} c_i^2(t),$$

(32)

where

$$c_i(t) = \int_0^1 f(x) \phi_i(x) dx = \int_0^1 f(x)(2i+1)^{\frac{1}{2}} P_i(2x-1) dx.$$

writing $2x - 1 = y$,

$$c_i(t) = \int_{-1}^{1} f\left(\frac{y+1}{2}\right)(2i+1)^{\frac{1}{2}} \frac{P_i(y)}{2} dy$$

$$= \left(\frac{2i+1}{2^2}\right)^{\frac{1}{2}} \int_{-1}^{1} f\left(\frac{y+1}{2}\right) P_i(y) dy$$

$$= \left(\frac{1}{2^2(2i+1)}\right)^{\frac{1}{2}} \int_{-1}^{1} f\left(\frac{y+1}{2}\right) d(P_{i+1}(y) - P_{i-1}(y)) dy$$

$$= -\frac{1}{4}\left(\frac{1}{(2i+1)}\right)^{\frac{1}{2}} \int_{-1}^{1} Df\left(\frac{y+1}{2}\right)(P_{i+1}(y) - P_{i-1}(y)) dy$$

$$= -\frac{1}{4}\left(\frac{1}{(2i+1)}\right)^{\frac{1}{2}} \int_{-1}^{1} Df\left(\frac{y+1}{2}\right) d\left(\frac{P_{i+2}(y) - P_i(y)}{2i+3} - \frac{P_i(y) - P_{i-2}(y)}{2i-1}\right) dy$$

$$= \frac{1}{8}\left(\frac{1}{(2i+1)}\right)^{\frac{1}{2}} \int_{-1}^{1} D^2 f\left(\frac{y+1}{2}\right)\left(\frac{P_{i+2}(y) - P_i(y)}{2i+3} - \frac{P_i(y) - P_{i-2}(y)}{2i-1}\right) dy.$$

Thus

$$|c_i(t)|^2 = \frac{1}{64}\left(\frac{1}{2i+1}\right)\left|\int_{-1}^{1} D^2 f\left(\frac{y+1}{2}\right)\left(\frac{(2i-1)P_{i+2}(y) - (4i+2)P_i(y) + (2i+3)P_{i-2}(y)}{(2i+3)(2i-1)}\right) dy\right|^2$$

$$\leq \frac{1}{64}\left(\frac{1}{2i+1}\right)\int_{-1}^{1}\left|D^2 f\left(\frac{y+1}{2}\right)\right|^2 dy \int_{-1}^{1}\left|\frac{(2i-1)P_{i+2}(y) - (4i+2)P_i(y) + (2i+3)P_{i-2}(y)}{(2i+3)(2i-1)}\right|^2 dy$$

$$\leq \frac{K^2}{32}\left(\frac{1}{2i+1}\right)\int_{-1}^{1}\left|\frac{(2i-1)P_{i+2}(y) - (4i+2)P_i(y) + (2i+3)P_{i-2}(y)}{(2i+3)(2i-1)}\right|^2 dy$$

$$< \frac{K^2}{32}\left(\frac{1}{2i+1}\right)\int_{-1}^{1} \frac{(2i-1)^2 P_{i+2}^2(y) + (4i+2)^2 P_i^2(y) + (2i+3)^2 P_{i-2}^2(y)}{(2i+3)^2(2i-1)^2} dy$$

$$= \frac{K^2}{32(2i+1)(2i-1)^2(2i+3)^2}\left[\frac{2(2i-1)^2}{2i+5} + \frac{8(2i+1)^2}{2i+1} + \frac{2(2i+3)^2}{2i-3}\right]$$

$$< \frac{K^2}{32(2i+1)(2i-1)^2(2i+3)^2} \frac{12(2i+3)^2}{(2i-3)}$$

$$< \frac{3K^2}{8(2i-3)^4}.$$

Therefore

$$\sum_{i=n}^{\infty} c_i^2(t) < \sum_{i=n}^{\infty} \frac{3K^2}{8(2i-3)^4}$$

$$= \frac{3K^2}{8} \sum_{i=n}^{\infty} \frac{1}{(2i-3)^4}$$

$$= \frac{3K^2}{8} \left(\frac{1}{96}\right) F_3\left(-\frac{3}{2}+n\right)$$

$$= \frac{K^2}{256} F_3\left(-\frac{3}{2}+n\right)$$

$$\|f(x) - S_n f(x)\|_{L^2}^2 \le \frac{K^2}{256} F_3\left(-\frac{3}{2}+n\right).$$

**Theorem:** If the constructed approximate solution $u_n$ satisfies

$$\left|D^2\left(\frac{\partial^\beta u_n}{\partial x^\beta}\right)\right|, \left|D^2\left(\frac{\partial^\gamma u_n}{\partial x^\gamma}\right)\right| \le K \quad \forall n \tag{33}$$

Then $u_n$ forms a weak asymptotic solution of the problem in Eq. (8) with initial condition in (9).

**Proof:** From the construction, it is clear that the sequence of approximate solutions

$$u_n(x,t) = \sum_{i=1}^{n} c_i(t)\phi_i(x), \text{ satisfy the following}$$

$$\frac{\partial^\alpha u_n(x,t)}{\partial t^\alpha} = -\nu S_n\left(\frac{\partial^\beta u_n(x,t)}{\partial x^\beta}\right) + k S_n\left(\frac{\partial^\gamma u_n(x,t)}{\partial x^\gamma}\right) + S_n(f(x,t)), \tag{34}$$

By substituting $u_n$ in (8), we get the following identity with the error term $\varepsilon_n(x,t)$.

$$\frac{\partial^\alpha u_n(x,t)}{\partial t^\alpha} = -\nu\left(\frac{\partial^\beta u_n(x,t)}{\partial x^\beta}\right) + k\left(\frac{\partial^\gamma u_n(x,t)}{\partial x^\gamma}\right) + S_n(f(x,t)) + \varepsilon_n(x,t) \tag{35}$$

Now subtracting Eq. (34) from Eq. (35) we get,

$$\varepsilon_n(x,t) = \nu\left(\left(\frac{\partial^\beta u_n}{\partial x^\beta}\right) - S_n\left(\frac{\partial^\beta u_n}{\partial x^\beta}\right)\right) + k\left(S_n\left(\frac{\partial^\gamma u_n}{\partial x^\gamma}\right) - \left(\frac{\partial^\gamma u_n}{\partial x^\gamma}\right)\right) \tag{36}$$

The error term in Eq.(36) can be estimated as,

$$\|\varepsilon_n(x,t)\|_{L^2} \leq \nu\left\|\left(\frac{\partial^\beta u_n}{\partial x^\beta}\right) - S_n\left(\frac{\partial^\beta u_n}{\partial x^\beta}\right)\right\|_{L^2} + k\left\|\left(\frac{\partial^\gamma u_n}{\partial x^\gamma}\right) - S_n\left(\frac{\partial^\gamma u_n}{\partial x^\gamma}\right)\right\|_{L^2} \tag{37}$$

Now assuming the following two conditions

$$\left|D^2\left(\frac{\partial^\beta u_n}{\partial x^\beta}\right)\right|, \left|D^2\left(\frac{\partial^\gamma u_n}{\partial x^\gamma}\right)\right| \leq K; \quad \forall n \tag{38}$$

Using Eq. (30) we get,

$$\left\|\left(\frac{\partial^\beta u_n}{\partial x^\beta}\right) - S_n\left(\frac{\partial^\beta u_n}{\partial x^\beta}\right)\right\|_{L^2}^2 \leq \frac{K^2}{256}F_3(-\frac{3}{2}+n), \quad \left\|\left(\frac{\partial^\gamma u_n}{\partial x^\gamma}\right) - S_n\left(\frac{\partial^\gamma u_n}{\partial x^\gamma}\right)\right\|_{L^2}^2 \leq \frac{K^2}{256}F_3(-\frac{3}{2}+n). \tag{39}$$

From Eq. (37) and Eq. (39), it follows

$$\|\varepsilon_n(x,t)\|_{L^2} \leq \frac{K}{16}\sqrt{F_3(-\frac{3}{2}+n)}(\nu+k). \tag{40}$$

So $\varepsilon_n(x,t) \to 0$ uniformly in the variables $x$ and $t$.

Multiplying test function $\phi$ on both sides of Eq. (35) and integrating with respect to $x$ in the interval (0, 1) we get

$$\int_0^1 \left(\frac{\partial^\alpha u_n(x,t)}{\partial t^\alpha} + \frac{\partial^\beta u_n(x,t)}{\partial x^\beta} - \lambda\frac{\partial^\gamma u_n(x,t)}{\partial x^\gamma} - f_n(x,t)\right)\phi(x)dx = \int_0^1 \varepsilon_n(x,t)\phi(x)dx, \tag{41}$$

Using Schwarz inequality in Eq. (41), we get,

$$\left|\int_0^1 \left(\frac{\partial^\alpha u_n(x,t)}{\partial t^\alpha} + \frac{\partial^\beta u_n(x,t)}{\partial x^\beta} - \lambda\frac{\partial^\gamma u_n(x,t)}{\partial x^\gamma} - f_n(x,t)\right)\phi(x)dx\right| \leq \|\varepsilon_n(x,t)\|\|\phi(x)\|, \tag{42}$$

Since $\|\varepsilon_n(x,t)\| \to 0$ as $n \to \infty$ uniformly. So $u_n$ forms a weak asymptotic solution.

## 5. Numerical results

**Example1:** Let us consider following advection-dispersion equation [30]

$$\frac{\partial^{\alpha} u(x,t)}{\partial t^{\alpha}} = -\frac{\partial^{\beta} u(x,t)}{\partial x^{\beta}} + \lambda \frac{\partial^{\gamma} u(x,t)}{\partial x^{\gamma}}, \qquad (43)$$

$t > 0,\ 0 < \alpha \leq 1,\ 0 < \beta \leq 1,\ 1 < \gamma \leq 2.$

with initial condition $u(x,0) = e^{-x}$.

Taking $\alpha = 1$, $\beta = 1$, $\gamma = 2$ and $\lambda = 1$, then Eq. (43) become

$$\frac{\partial u(x,t)}{\partial t} = -\frac{\partial u(x,t)}{\partial x} + \frac{\partial^2 u(x,t)}{\partial x^2}, \qquad (44)$$

whose exact solution is given by $u(x,t) = e^{-x+2t}$. In Fig. 1.1 $E1, E2, E3, E4, E5$ represent absolute error corresponding to the time $t = 0.00001, 0.1, 0.5, 0.9,\ 0.99999$ respectively.

In Fig. 1.2, 1.3 and 1.4 we can see that as the value of n (no. of basis element) increases, the approximate solutions come closer and closer to the exact solution at different time t=0.5, 0.8 and 1 respectively.

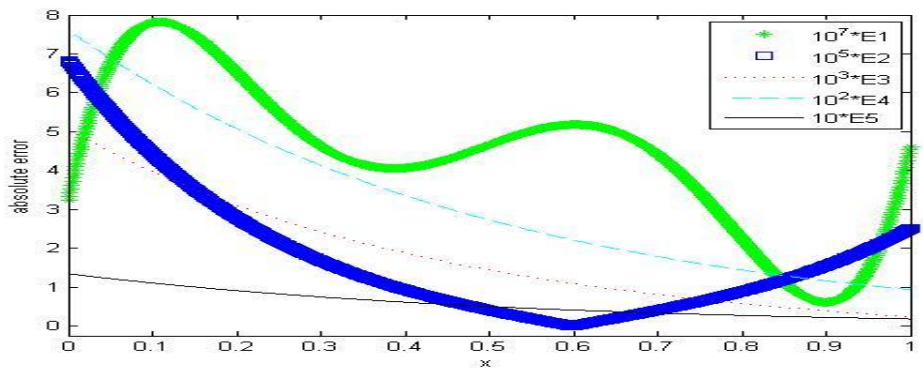

**Fig. 1.1 Comparison of the absolute errors at different time level**

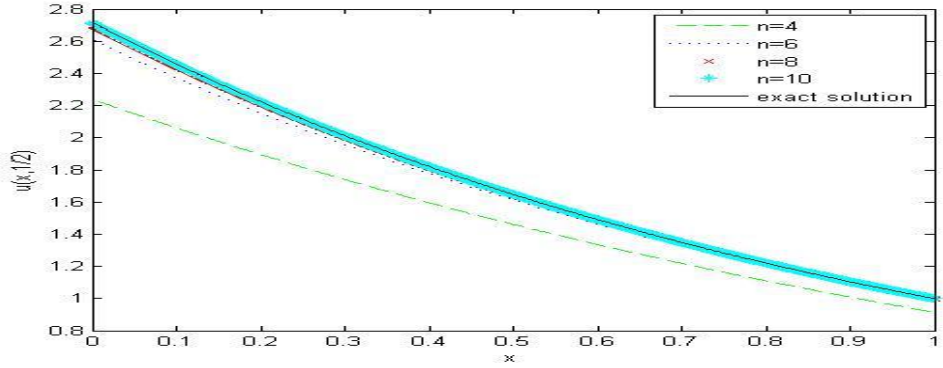

**Fig. 1.2** Comparison of numerical and exact solution at different values of n at t=0.5.

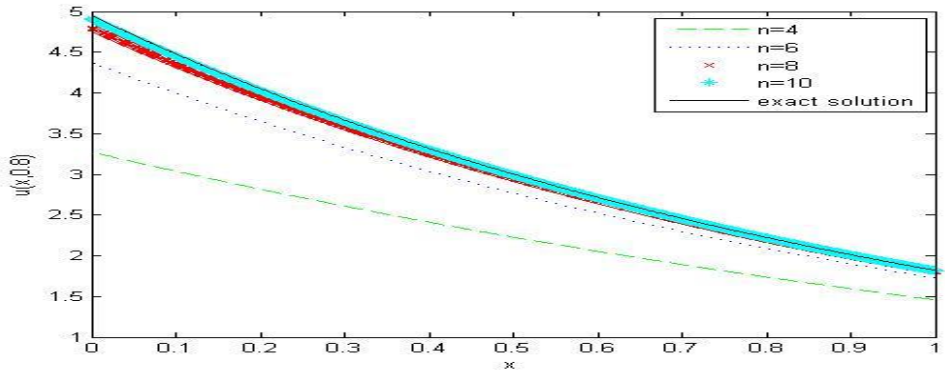

**Fig. 1.3** Comparison of numerical and exact solution at different values of n at t=0.8.

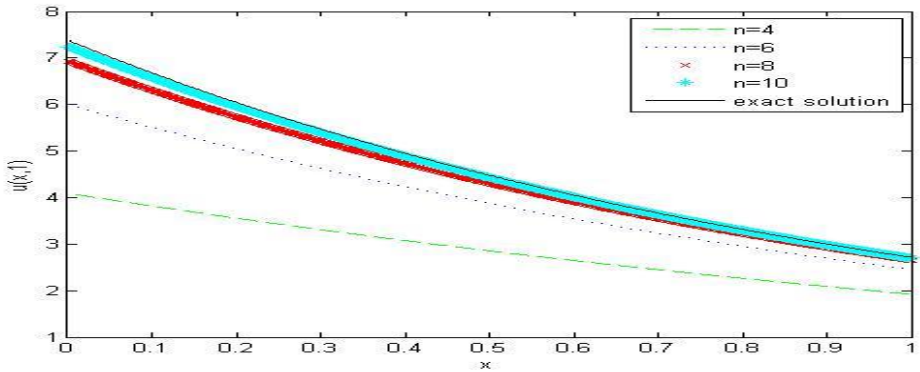

**Fig. 1.4** Comparison of numerical and exact solution at different values of n at t=1.

**Example2:** We consider the following non-homogeneous fractional advection-dispersion equation [30]

$$\frac{\partial^{2\beta} u(x,t)}{\partial x^{2\beta}} - \frac{\partial u(x,t)}{\partial x} = \frac{\partial u(x,t)}{\partial t} + f(x,t), \qquad (45)$$

$t > 0,\ x > 1,\ 0 < \beta \leq 1.$ With initial condition $u(x,0) = x^2,\ f(x,t) = 2 - 2t - 2x$ and $\beta = 1$, Eq. (45) has exact solution $u(x,t) = x^2 + t^2$.

We have also discussed problem Eq. (45) with different values of $\beta = 0.125, 0.25, 0.375, 0.5, 1.$ and different non-homogeneous terms. Fig. 2.1 and 2.2; $E1$, $E2$, $E3$, $E4$, $E5$ show absolute error for $\beta = 0.125, 0.25, 0.375, 0.5, 1$ respectively at the time levels t=0.1 and 0.5.

In Fig. 2.3, 2.4 and 2.5 $E1, E2, E3, E4, E5$ show absolute error for $t = 0.1, 0.25, 0.5, 0.75, 0.9$ respectively at fix values of $\beta = 0.25, 0.375, 1$.

From the graph it is clear that as the value of $t$ decreases error decreases, when it goes near to zero error is in $10^{-7}$ level.

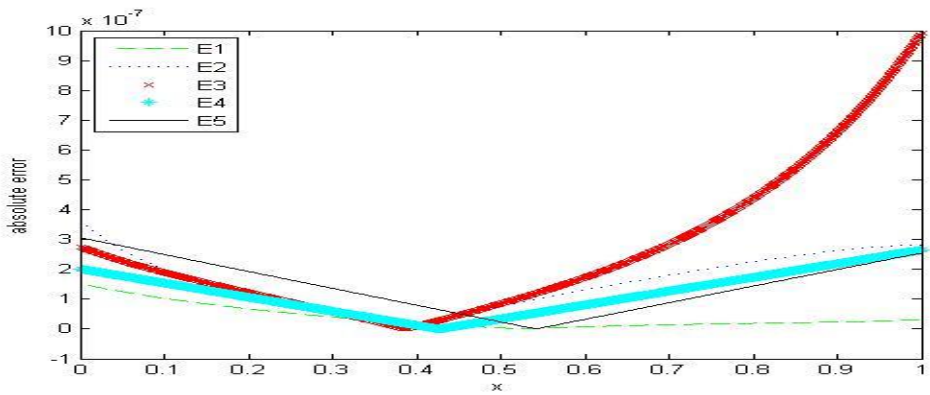

**Fig. 2.1 Comparison of the absolute errors at different values of beta and t=0.1.**

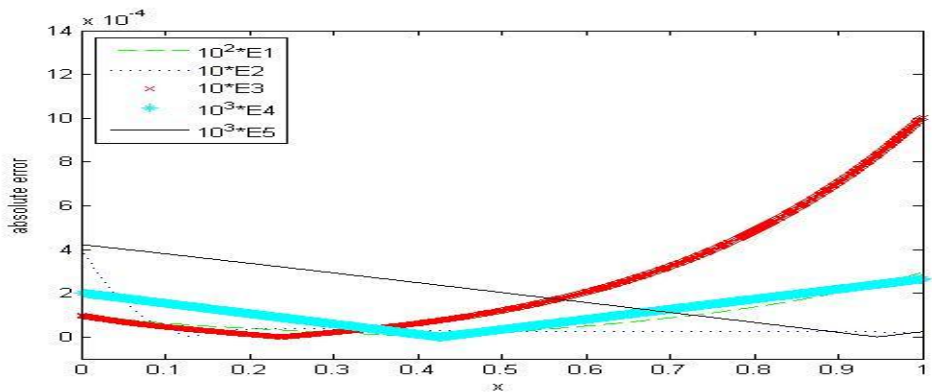

**Fig. 2.2 Comparison of the absolute errors at different values of beta and t=0.5.**

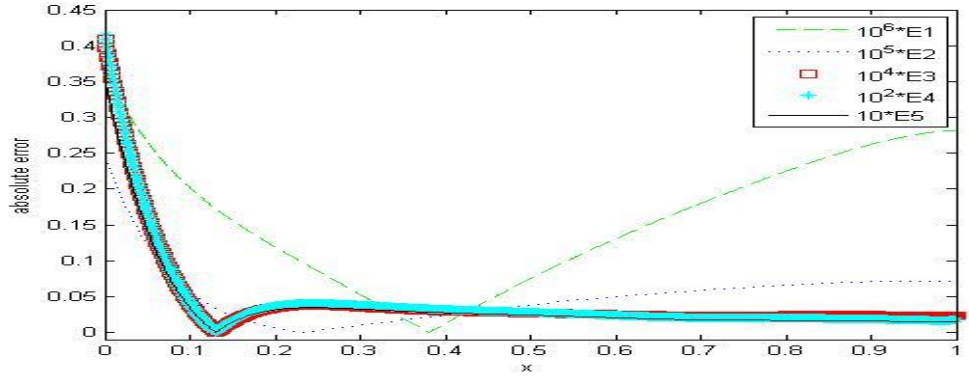

**Fig. 2.3 Comparison of the absolute errors at different values of time and beta=0.25.**

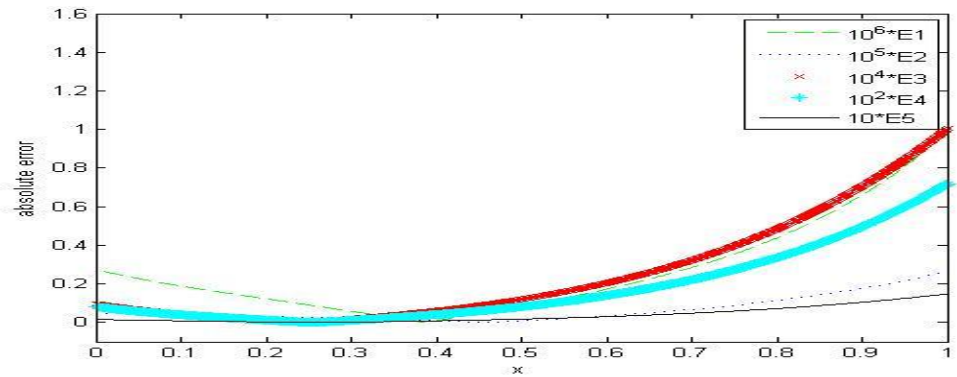

**Fig. 2.4 Comparison of the absolute errors at different values of time and beta=0.375.**

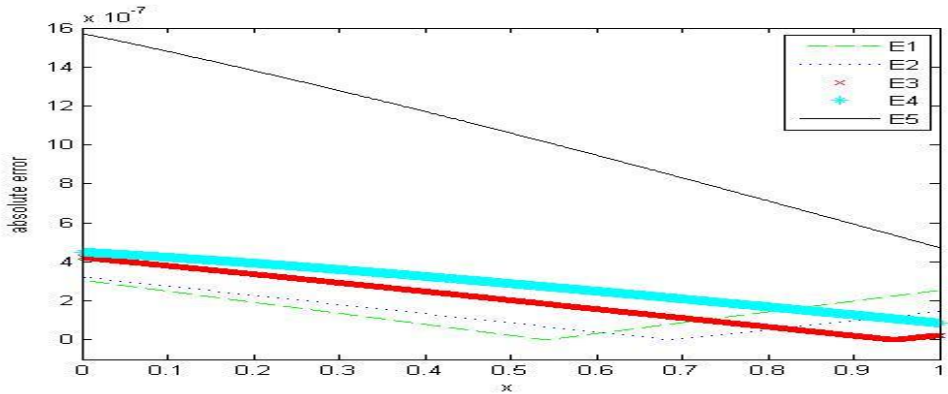

**Fig. 2.5 Comparison of the absolute errors at different values of time and beta=1.**

**Example3:** We consider the following non-homogeneous time-fractional advection-dispersion equation

$$\frac{\partial^\alpha u(x,t)}{\partial t^\alpha} = \frac{\partial^2 u(x,t)}{\partial x^2} - \frac{\partial u(x,t)}{\partial x} + f(x,t), \qquad (46)$$

$t > 0, x > 0, 0 < \alpha \leq 1.$

With initial condition $u(x,0) = x^2$, $f(x,t) = 2x - 2 + \frac{\sqrt{2}}{\sqrt{\frac{3}{2}}} t^{\frac{1}{2}}$ and $\alpha = \frac{1}{2}$. Eq.(46) has exact solution $u(x,t) = x^2 + t$.

We also solved the above equation for different values of $\alpha = 0.25, 0.5, 0.75, 1$ and different non-homogeneous terms. In Fig. 3.1, 3.2, 3.3 and 3.4; $E1, E2, E3, E4, E5$ represent absolute error for $t = 0.1, 0.25, 0.5, 0.75, 0.9$ respectively at fixed values of $\alpha = 0.25, 0.5, 0.75$ and $1$.

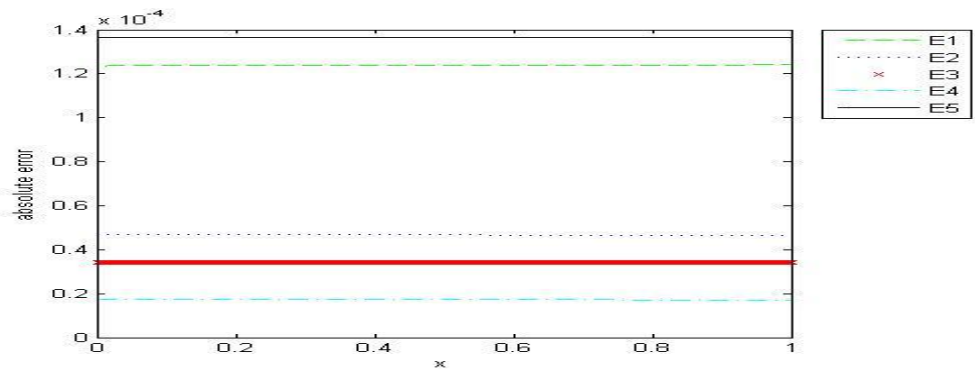

**Fig. 3.1 Comparison of the absolute errors at different values of time and alpha=0.25.**

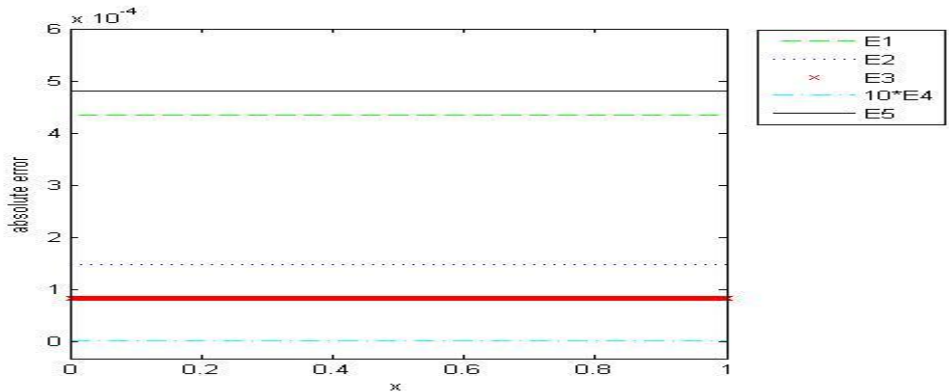

**Fig. 3.2 Comparison of the absolute errors at different values of time and alpha=0.5.**

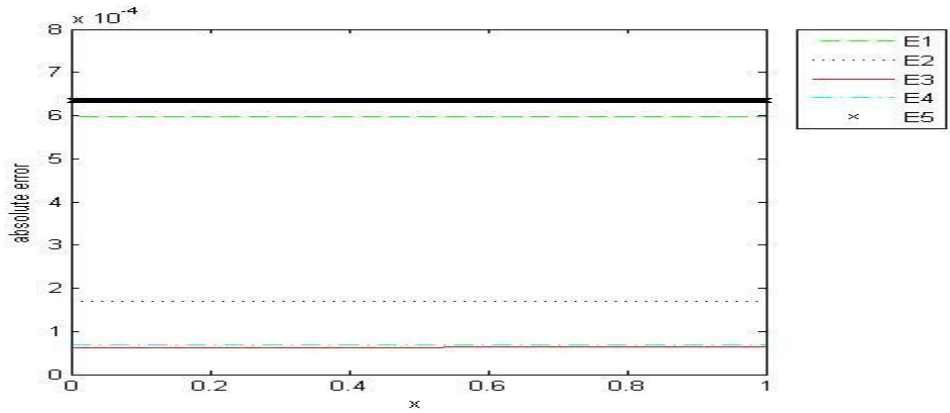

**Fig. 3.3 Comparison of the absolute errors at different values of time and alpha=0.75.**

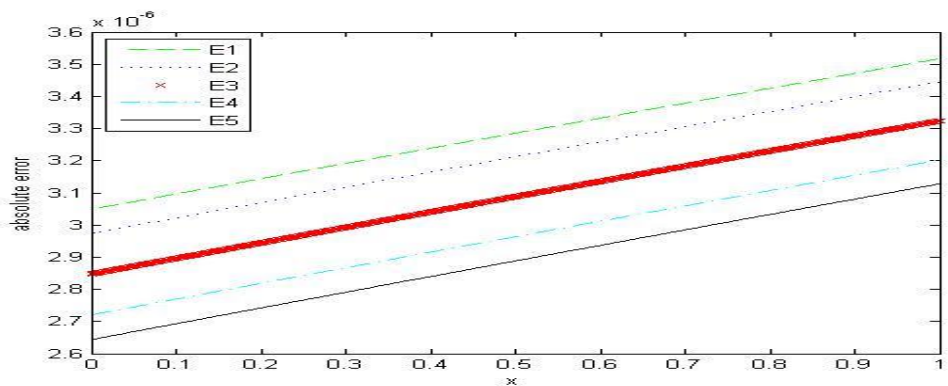

**Fig. 3.4 Comparison of the absolute errors at different values of time and alpha=1.**

**Conclusions:** In this paper we used Galerkin type approximation to construct approximate solution to fractional advection-dispersion equation. We showed convergence of solution assuming some conditions on the approximate solutions. Numerically those conditions are ideal conditions. Our algorithm shows good convergence rate to the exact solution. We carried out our method in the domain $[0,1] \times [0,\infty]$. This method can be used in any bounded domain in x-axis, namely $[-k,k] \times [0,\infty]$, by scaling Legendre functions carefully.